\input Tex-document.sty

\def\frac#1#2{{#1\over #2}}
\let\phi=\varphi\def\<{\langle}\def\>{\rangle}
\def\A{{\cal A}}\def\B{{\cal B}}
\def\C{{\bf C}}\def\epsilon{\varepsilon}
\def\H{{\cal H}}
\def\L{{\cal L}}
\def\M{{\cal M}}\def\N{{\bf N}}

\def\R{{\bf R}}
\def\TTT{{\cal T}}
\def\Z{{\bf Z}}

\def\II{{\rm II}} \font\eightrm=cmr8 \font\eightsl=cmsl8
\font\eightbf=cmbx8  
\def\tensor{\mathop{\bar\otimes}}  
\def\iint{\int\!\!\!\int}

\def\firstpage{787}

\pageno=787

\def\shorttitle{Free Probability, Free Entropy and Applications to von Neumann Algebras}
\def\shortauthor{Liming Ge}

\title{\centerline{Free Probability, Free Entropy and}
\centerline{Applications to von Neumann Algebras}}

\author {Liming Ge\footnote{\eightrm *}{\eightrm Academy of Mathematics and System Science, CAS, Beijing 100080,
China. Department of Mathematics UNH, Durham, NH 03824, USA.
E-mail: liming@math.unh.edu}}

\vskip 7mm



 This talk is organized as follows: First we explain some basic concepts in
non-commutative probability theory in the frame of operator algebras. In Section 2, we discuss related topics in
von Neumann algebras. Sections 3 and 4 contain some of the key ideas and results in free probability theory. Last
section states some of the important applications of free probability theory.





\head{1. Non-commutative probability spaces}

In general, a non-commutative probability space is a pair $(\A, \tau)$, where $\A$ is a unital algebra (over the
field of complex numbers $\C$) and $\tau$ a linear functional with $\tau(I)=1$, where $I$ is the identity of $\A$.
Elements of $\A$ are called {\it random variables}. Since positivity is a key concept in (classical) probability
theory, this can be captured by assuming that $\A$ is a * algebra and $\tau$ is positive (i.e., a state). Elements
of the form $A^*A$ are called {\it positive} (random variables).

A state $\tau$ is a {\it trace} if $\tau(AB)=\tau(BA)$. We often
require that $\tau$ be a faithful trace ($\tau$ corresponds to the
classical probability measure, or the integral given by the
measure). In this talk, we always assume that $\A$ is a unital *
algebra over $\C$ and $\tau$ a faithful state on $\A$. Subalgebras
of $\A$ are always assumed unital * subalgebras.

Examples of noncommutative probability spaces often come from
operator algebras on a Hilbert space and the states used here are
usually vector states.

A {\it C*-probability space} is a pair $(\A, \tau)$, where $\A$ is
a unital C*-algebra (norm closed subalgebra of $\B(\H)$) and
$\tau$ is a state on $\A$.

A {\it W*-probability space} is a pair $(\M,\tau)$ consisting of a
von Neumann algebra $\M$ (strong-operator closed C*-subalgebra of
$\B(\H)$) and a {\it normal\/} (i.e., countably additive) state
$\tau$ on $\M$.

The following are some more basic concepts:

\noindent {\bf  Independence:} In a noncommutative probability
space $(\A, \tau)$, a family $\{\A_j\}$ of subalgebras $\A_j$ of
$\A$ is {\it independent} if the subalgebras commute with each
other and, for $n\in\N$, $\tau(A_1\cdots
A_n)=\tau(A_1)\cdots\tau(A_n)$ for all $A_k$ in $\A_{j_k}$ and
$j_k\neq j_l$ whenever $k\neq l$.

This independence gives a ``tensor-product'' relation among
subalgebras $\A_j$: if $\A$ is generated by $\A_j$, then
$\A\cong\otimes_j\A_j$ (in the case of C*- or W*-probability
spaces, the tensor-product shall reflect the corresponding
topological structures on $\A$ and $\A_j$).

\noindent {\bf  Distributions and moments:} Given $(\A, \tau)$,
for $A$ in $\A$, we define a map $\mu_A: \C[x]\to \C$ by
$\mu_A(p(x))=\tau(p(A))$. Then $\mu_A$ is the {\it distribution}
of $A$. For $A_1,\ldots,A_n$ in $\A$, the {\it joint distribution}
$\mu_{A_1,\ldots,A_n}: \C\<x_1,\ldots,x_n\>\to\C$ is given by $$
\mu_{A_1,\ldots,A_n}(p(x_1,\ldots,x_n))=\tau(p(A_1,\ldots,A_n)).
$$ If $p$ is a monomial, $\tau(p(A_1,\ldots,A_n))$ is called a
{\it (p-)moment}. When random variables are non self-adjoint, one
also considers (joint) * distributions of random variables, that
can be defined in a similar way. In this case, there is a natural
identification of $\C\<x_1,\ldots,x_n,x_1^*,\ldots,x_n^*\>$ with
the semigroup algebra $\C S_{2n}$, where $S_{2n}$ is the free
semigroup on $2n$ generators. Monomials are given by words in
$S_{2n}$.

\noindent {\it Conditional Expectations:} Suppose $\B$ is a
subalgebra of $\A$. A {\it conditional expectation} from $\A$ onto
$\B$ is a $\B$-bimodule map (a projection of norm one in the case
of C*-algebras) of $\A$ onto $\B$ so that the restriction on $\B$
is the identity map.

Many other concepts in probability theory and measure theory can
be generalized to operator algebras, especially von Neumann
algebras which can be regarded as non-commutative measure spaces.
For basic operator algebra theory, we refer to [KR] and [T].

\head{2. GNS representation and von Neumann algebras}

Given a C*-probability space $(\A, \tau)$, one defines an inner product $\<A,B\>=\tau(B^*A)$ on $\A$. Let
$L^2(\A,\tau)$ be the Hilbert space obtained by the completion of $\A$ under the $L^2$-norm given by this inner
product. Then $\A$ acts on $L^2(\A,\tau)$ by left multiplication. This representation of $\A$ on the Hilbert space
$L^2(\A,\tau)$ is called the {\it GNS representation}. In a similar way, one can define $L^p(\A,\tau)$, where
$\|A\|_p=\tau(|A|^p)^{1/p} =\tau((A^*A)^{p/2})^{1/p}$. The von Neumann algebra generated by $\A$ (or the
strong-operator closure of $\A$) is sometimes denoted by $L^\infty(\A,\tau)$ ($\subset L^p(\A,\tau)$, $p\ge 1$).
All von Neumann algebras admit such a form. Any von Neumann algebra is a (possibly, continuous) direct sum of
``simple'' algebras, or factors (algebras with a trivial center). Von Neumann algebras that admit a faithful
(finite) trace are said to be {\it finite\/}. The classification of (infinite-dimensional) finite factors has
become the central problem in von Neumann algebras.

Murray and von Neumann [MN] also separate factors into three
types:

Type I: Factors contain a minimal projection. They are isomorphic
to full matrix algebras $M_n(\C)$ or $\B(\H)$.

Type II: Factors contain a ``finite'' projection but without
minimal projections: it is said to be of type $\II_1$ when the
identity $I$ is a finite projection; of type $II_\infty$ when $I$
is infinite. Every type $\II_\infty$ is the tensor product of
$\B(\H)$ with a factor of type $\II_1$.

Type III: Every (non zero) projection is infinite.

\noindent {\bf Examples of von Neumann algebras:} 1) Let $\A=\C G$
for some discrete group $G$, $\H$ be $l^2(G)$ and $\tau$ be the
vector state given by the vector that takes value 1 at $g$ and 0
elsewhere. Then $\tau$ is a trace, $l^2(G)=L^2(\A,\tau)$ and the
weak (or strong) operator closure of $\A$ is called the group von
Neumann algebra, denoted by $\L_G$. We have that $\L_G$ is a
factor if and only if each conjugacy class of $G$ other than the
identity is infinite (i.c.c.). For example, the free group $F_n$
($n\ge2$, on $n$ generators) is such an i.c.c. group.

2) Suppose $(\Omega,\mu)$ is a measure space with a
$\sigma$-finite measure $\mu$, $G$ is a group and $\alpha$ is a
measurability preserving action of $G$ on $\Omega$. Formally, we
have an algebra $L^\infty(\Omega,\mu)G$ ($=\A$) similar to the
group algebra definition:
$(\phi(x)g)(\psi(x)h)=\phi(x)\psi(g^{-1}(x))gh$, for $\phi,\psi\in
L^\infty(\Omega,\mu)$ and $g,h\in G$. Assume that $G$ acts freely
(i.e., for any $g$ in $G$ with $g\neq e$, the set $\{x\in\Omega:
g(x)=x\}$ has measure zero). Define an action of $\A$ on the
Hilbert space $\bigoplus_{g\in G}L^2(\Omega,\mu)g$ by left
multiplication ( which is induced by the multiplication in $\A$),
where $L^2(\Omega,\mu)g$ is an isomorphic copy of
$L^2(\Omega,\mu)$. Then the von Neumann algebra generated by $\A$
is called the cross product von Neumann algebra, denoted by
$L^\infty(\Omega,\mu)\times_\alpha G$. This cross product is a
factor if and only if $\alpha$ is ergodic. If $\alpha$ is ergodic
and {\it not} a measure preserving action, then
$L^\infty(\Omega,\mu)\times_\alpha G$ is a factor of type III.
Type II factors are obtained from measure preserving actions (with
$\Omega$ a non atomic measure space) and the finiteness of $\mu$
gives rise to type $\II_1$ factors.

The above 1) and 2) are the two basic constructions of von Neumann
algebras (given by Murray and von Neumann [MN]). A. Connes [C]
shows that there are finite factors that cannot be constructed by
1). It was a longstanding open problem whether every (finite) von
Neumann algebra can be obtained by using the construction in 2).
Using free probability theory, especially the notion of free
entropy, Voiculescu [V2] gives a negative answer to this question.
We shall discuss some details later in the talk.

In recent years, the focus of studies of von Neumann algebras is
centered on factors of type $\II_1$. Many of the unsolved problems
in operator algebras are also reduced to this class. We end this
section with two of the (still) open problems from the list of 20
questions asked by Kadison in 1967.

\noindent{\bf 1.} {The weak-operator closure of the left regular
representation of the free (non-abelian) group on two or more
generators is a factor of type $\II_1$. Are these factors
isomorphic for different numbers of generators?}

\noindent{\bf 2.} {Is each factor generated by two self-adjoint
operators? ---each von Neumann algebra? ---the factor arising
from the free group on three generators? ---is each von Neumann
algebra finitely generated?}

Three of those 20 problems were answered in the last ten years by
using free probability and free entropy. We explain some of the
theory involved in the following two sections.

\head{3. Free independence}

Suppose $(\A,\tau)$ is a C*-probability space. We assume that $\tau$ is a trace. A family $\A_\iota$, $\iota\in\bf
I$, of unital
* subalgebras of $\A$ are called {\it free} with respect to $\tau$ if
$\tau(A_{1}A_{2}\cdots A_{n})=0$ whenever $A_{j}\in \A_{\iota_j}$,
$\iota_1\neq\cdots\neq\iota_n$ ($\iota_1$ and $\iota_3$ may be the
same) and $\tau(A_j)=0$ for $1\le j\le n$ and every $n$ in $\N$. A
family of subsets (or elements) of $\A$ are said to be {\it free}
if the unital * subalgebras they generate are free.

Note that freeness is a highly noncommutative notion, the
non-commutativity (or algebraic freeness) of free random variables
is encoded in the definition. Recall some basic concepts in free
probability theory.

\noindent {\bf Semicircular elements:} The Gaussian laws in
classical theory is replaced by the semicircular laws. The {\it
semicircular law} centered at $a$ and of radius $r$ is the
distribution $\mu_{a,r}: \C[x]\to \C$ such that $$
\mu_{a,r}(\phi(x))=\frac 2{\pi r^2}\int_{a-r}^{a+r}\phi(t)\sqrt{
r^2-(t-a)^2}dt, $$ for each $\phi(x)$ in $\C[x]$. A self-adjoint
random variable $A$ in $(\A,\tau)$ is said to be {\it (standard)
semicircular} if its distribution is $\mu_{0,1}$. An element
$X=A+iB$ is {\it circular} if $A$ and $B$ are free semicircular.
The following theorem is proved by D. Voiculescu (see [VDN]).

\noindent {\bf Free Central Limit Theorem:} {\it Let
$\{A_j\}_{j=1}^\infty$ be a free family of identically distributed
random variables in $(\A,\tau)$ with $\tau(A_j)=0$ and
$\tau(A_j^2)=\frac{r^2}4$ for some positive number $r$. Then the
distribution of $\frac{A_1+\cdots+A_n}{\sqrt n}$ converges to the
semicircular law $\mu_{0,r}$.} \vskip4pt

Classical independence corresponds to tensor products; while the
above free independence introduced by Voiculescu is given by
certain free products. Recall some examples of such freeness.

\noindent {\bf Free products:} If $G=G_1*G_2$, then $\L_{G_1}$
and $\L_{G_2}$ are free in $\L_G$, here the trace is the one
given by the unit vector associated with any group element.

Let $(\A_1,\tau_1)$ and $(\A_2, \tau_2)$ be two W*-probability
spaces. Suppose $\A_0$ is the (amalgamated algebraic) free product
of $\A_1$ and $\A_2$ (over the scalars). Then there is a unique
$\tau$ on $\A_0$ such that $\A_1$ and $\A_2$ are free with respect
to $\tau$ and the restrictions of $\tau$ on $\A_1$ and $\A_2$
equal to $\tau_1$ and $\tau_2$, respectively. Let $\A$ be the weak
operator closure of $\A_0$ acting on $L^2(\A_0,\tau)$. Then $\A$
is called the (reduced von Neumann algebra) free product of $\A_1$
and $\A_2$ (with respect to $\tau$), denoted by $\A_1*\A_2$ (and
$\tau=\tau_1*\tau_2$). For example, $\L_{F_2}\cong L^\infty[0,1]*
L^\infty[0,1]$.

\noindent {\bf Full Fock space construction:} Let $\H_0$ be a real
Hilbert space and $\H$ be $\H_0\otimes\C$. Its {\it full Fock
space} is $$ \TTT(\H)=\C 1\oplus\bigoplus_{n\ge1}\H^{\otimes n}.
$$ For $h\in\H_0$, let the {\it left creation operator}
$l(h)\in\B(\TTT(\H))$ be given by $l(h)\xi=h\otimes\xi$. Then
$l(h)^*1=0$ and $l(h)^*\xi_1\otimes \xi_2\otimes\cdots\otimes
\xi_n=\<\xi_1, h\> \xi_2\otimes \cdots \otimes \xi_n$, so
$l(h_1)^*l(h_2)=\<h_2,h_1\>I$. Let $C^*(l(\H_0))$ (or
$W^*(l(\H_0))$) be the C*- (or W*-) algebra generated by $\{l(h)|
h\in\H_0\}$. Let $\omega_\H$ be the vector state given by the
vector $1\in\TTT(\H)$. Here $1$ is called the {\it vacuum vector}
and $\omega_\H$ the {\it vacuum state}. If $\H_1$ and $\H_2$ are
orthogonal subspaces of $\H_0$, then $C^*(l(\H_1))$ and
$C^*(l(\H_2))$ are free (with respect to $\omega_\H$). If $h$ is a
unit vector in $\H_0$, then $(l(h)+l(h)^*)/2$ is semicircular with
distribution $\mu_{0,1}$. \vskip4pt

\noindent {\bf Gaussian Random Matrices:} Let
$X(s,n)=(f_{ij}(s,n))$ in $M_n(L^\infty[0,1])$ be real random
matrices, where $n\in\N$ and $s\in S$ for some index set $S$.
Assume that $f_{ij}(s,n)=f_{ji}(s,n)$ and $\{f_{ij}(s,n): i,j,s\}$
(given each $n$) is a family of independent Gaussian $(0,1/n)$
random variables. Let $D_n$ be a constant diagonal matrix in
$M_n(\R)$ having a limit distribution (as $n\to\infty$). Then
$\{X(s,n)\}\cup\{D_n\}$ is asymptotically free as $n\to\infty$ and
$\{X(s,n): s\in S\}$ converges in distribution to a free
semicircular family.  As a corollary, Voiculescu shows that
$\L_{F_5}\otimes M_2(\C)\cong \L_{F_2}$. Moreover $\L_{F_r}\otimes
M_n(\C)\cong\L_{F_{1+\frac{r-1}{n^2}}}$, for any real number $r$,
and $\L_{F_r}*\L_{F_s}\cong\L_{F_{r+s}}$. Now we know that either
$\L_{F_r}$, $r>1$, are all isomorphic to each other or they are
all non isomorphic factors (see [D] and [R]).

Further studies of free probability theory have been pursued by
many people in several directions, such as infinitely divisible
laws, free brownian motion, etc. (we refer to [B], [VDN] and [HP]
for details).

\head{4. Free entropy}

Free entropy is a non commutative analogue of classical Shannon entropy. First we recall the definition of entropy
and its basic properties.

\noindent {\it Classical Entropy:} Let $(\Omega, \Sigma, \mu)$ be
a probability space with probability measure $\mu$ and
$f_1,\ldots, f_n:\Omega\to\R$ be random variables. Suppose $\phi$
is the density function on $\R^n$ corresponding to the joint
distribution of $f_1,\ldots,f_n$. Then the entropy: $$
H(f_1,\ldots,f_n)=-\int_{\R^n}
\phi(t_1,\ldots,t_n)\log\phi(t_1,\ldots, t_n)dt_1\cdots dt_n. $$
Here are two important properties of entropy:
$H(f_1,\ldots,f_n)=H(f_1)+\cdots+H(f_n)$ if and only if
$f_1,\ldots,f_n$ are independent; when assuming that $E(f_j^2)=1$,
$H$ is maximal if and only if $f_1,\ldots,f_n$ are Gaussian
independent $(0,1)$ random variables.

\noindent {\it Free entropy:} Let $X_1,\ldots, X_n$ be
self-adjoint random variables in $(\A,\tau)$. For any
$\epsilon>0$, when $k$ large, there may be self-adjoint matrices
$A_1,\ldots, A_n$ in $M_k(\C)$ such that {\it ``the algebra
generated by $X_j$'s looks like the algebra generated by $A_j$'s
within $\epsilon$.''} More precisely, for any $\epsilon>0$, large
$m\in\N$ and any monomial $p$ in $\C\<x_1,\ldots,x_n\>$ with
degree less than or equal to $m$, choose $k$ large enough so that
$$
|\tau(p(X_1,\ldots,X_n))-\tau_k(p(A_1,\ldots,A_n))|<\epsilon.\leqno(*)
$$

Let $\Gamma_R(X_1,\ldots,X_n; m,k,\epsilon)$ be the set of all
self-adjoint matrices $(A_1,\ldots,A_n)$ in $M_k(\C)$, with
$\|A_j\|\le R$, such that $(*)$ holds.
 The limit of the ``normalized''
measurement of $\Gamma_R(X_1,\ldots,X_n; m,k,\epsilon)$ is called the free entropy of $X_1,\ldots,X_n$. Voiculescu
[V1] shows that this limit is independent of $R$ when it is larger than $\max\{\|X_j\|: j=1,\ldots,n\}$. Here we
will fix such a constant $R$ and use $\Gamma(X_1,\ldots,X_n$; $m,k,\epsilon)$ to denote $\Gamma_R(X_1,\ldots,X_n;
m,k,\epsilon)$. Let $vol$ be the euclidean volume in real euclidean space $(M_k(\C)^{s.a.})^n$ (here ``s.a.''
denote the self-adjoint part and the euclidean norm $\|A\|_e^2={\rm Tr}(A^2)$). Now we define, successively,
$$ \eqalign{ \chi(X_1,\ldots, X_n; m,k,\epsilon)
                       &=\log vol(\Gamma(X_1, \ldots,
               X_n; m, k,\epsilon)),\cr
\chi(X_1,\ldots, X_n; m,\epsilon)
                       &=\limsup_{k\to\infty}(k^{-2}\chi(X_1,
            \ldots,X_n; m,k,\epsilon)       +\frac n2\log k),\cr
\chi(X_1,\ldots,X_n)&=\inf \{\chi(X_1,\ldots, X_n; m,\epsilon) :
m\in\N, \epsilon>0\}.\cr } $$ We call $\chi(X_1,\ldots,X_n)$ the
{\it free entropy\/} of $(X_1,\ldots,X_n)$.

The following are some basic properties of free entropy (proved in
[V1]):

{\rm (i)} $\chi(X_1,\ldots, X_n)\le \frac n2\log(2\pi
eC^2n^{-1})$;

{\rm (ii)} $\chi(X_1)=\iint \log|s-t|d\mu_1(s)d\mu_1(t)+\frac34
+\frac12\log2\pi$, where $\mu_1$ is the (measure on the spectrum
of $X_1$ corresponding to the) distribution of $X_1$;

{\rm (iii)} $\chi(X_1,\ldots, X_n)=\chi(X_1)+\cdots+ \chi(X_n)$
when $X_1, \ldots, X_n$ are free random variables.

Voiculescu also introduces a notion of free entropy dimension
$\delta(X_1,\ldots,X_n)$ which is given by
$$\delta( X_1,\ldots,X_n)
=n+\limsup_{\epsilon\to0}\frac {\chi(X_1+\epsilon S_1,\ldots,
X_n+\epsilon S_n: S_1,\ldots,S_n)}{|\log\epsilon|}
$$ where $\{S_1,\ldots,S_n\}$ is a standard semicircular family and
$\{X_1,\ldots,X_n\}$ and $\{S_1, \ldots$, $S_n\}$ are free.

\head{5. Applications}

In this section, we review some of the applications of free entropy in von Neumann algebras.

\noindent{\bf Theorem 1.} ([V1]) \hskip8pt{\it Let
$S_1,\ldots,S_n$ be a standard (0,1) free semicircular family.
Then}
 $$\eqalign{
 \chi(S_1,\ldots,S_n)& =n(\frac12+\log\frac{\sqrt{2\pi}}2);\cr
 \delta(S_1,\ldots,S_n)& =n.\cr}
 $$

Note that $S_1,\ldots,S_n$ generate $\L_{F_n}$ as a von Neumann
algebra.

\noindent{\bf Theorem 2.} ([V2])\hskip8pt {\it If a finite von
Neumann algebra $\M$ (with a trace) has a Cartan subalgebra, then,
for any self-adjoint generators $X_1,\ldots,X_n$ of $\M$,
$\delta(X_1,\ldots,X_n)\le1$.}

This implies that free group factors do not contain any Cartan
subalgebras which answers a longstanding question in the subject.

The following result generalizes Voiculescu's result and answers
Problem 11 on Kadison's 1967 problem list (unpublished).

\noindent{\bf Theorem 3.} ([G1])\hskip8pt {\it If a finite von
Neumann algebra $\M$ (with a trace) has a simple maximal abelian
subalgebra, then, for any self-adjoint generators $X_1,\ldots,X_n$
of $\M$, $\delta(X_1,\ldots,X_n)\le2$.}

Furthermore, we prove the following result which shows the
existence of a separable prime factor (the one that is not the
tensor product of two factors of the same type).

\noindent{\bf Theorem 4.} ([G2]) {\it If $\M=\M_1\tensor\M_2$ for
some infinite dimensional finite von Neumann algebras $\M_1$ and
$\M_2$, then, for any self-adjoint generators $X_1,\ldots,X_n$ of
$\M$, $\delta(X_1,\ldots,X_n)\le1$.}

It is an outstanding open question whether the free entropy
dimension is an invariant for von Neumann algebras. Following
Voiculescu's result in [V3], we prove the following:

\noindent{\bf Theorem 5.} ([GS])\hskip4pt {\it For any
self-adjoint generators $X_1,\ldots,X_n$ of $\L_{SL_m(\Z)}$,
$m\ge3$, $\delta(X_1,\ldots,X_n)\le1$.}

The above theorem is not true for $SL_2(\Z)$.

\parindent 12mm

\head{References}

\item{[B]} P. Biane, Proceedings of ICM2002.

\item{[C]} A. Connes, {\eightsl  Sur la classification des facteurs
de type $\II_1$}, (French. English summary) C. R. Acad. Sci. Paris
S\'er. A-B {\eightbf 281} (1975), A13--A15.

\item{[D]} K.\ Dykema, {\eightsl  Free products of hyperfinite von
Neumann algebras and free dimension,} Duke Math. J., {\eightbf 69}
(1993), 97--119.

\item{[G1]} L.\ Ge, {\eightsl Applications of free entropy to finite
von Neumann algebras,} Amer.\ J.\ Math., {\eightbf 119} (1997),
467--485.

\item{[G2]} L.\ Ge, {\eightsl Applications of free entropy to finite
von Neumann algebras, II,} Ann.\ of Math., {\eightbf 147} (1998),
143--157.

\item{[GS]} L.\ Ge and J.\ Shen, {\eightsl Free entropy and property
T factors}, Proc.\ Nat.\ Acad.\ Sci.\ (USA), {\eightbf 97} (2000),
9881--9885.

\item{[HP]} F. Hiai and D. Petz,  {\eightsl ``The Semicircle Law,
Free Random Variables and Entropy''}, Mathematical Surveys and
Monographs, {\eightbf 77}, American Mathematical Society,
Providence, RI, 2000.

\item{[KR]} R.\ Kadison and J.\ Ringrose, {\eightsl ``Fundamentals of
the Operator Algebras,''} vols. I and II, Academic Press, Orlando,
1983 and 1986.

\item{[MN]} F. J. Murray and J. von Neumann, {\eightsl On rings of
operators}, Ann. of Math., {\eightbf 37} (1936), 116--229.

\item{[R]} F.\ Radulescu, {\eightsl Random matrices, amalgamated free
products and subfactors of the von Neumann algebra of a free
group, of noninteger index}, Invent. Math., {\eightbf 115} (1994),
347--389.

\item{[T]} M. Takesaki, {\eightsl ``Theory of Operator Algebras}, vol
I, Springer-Verlag, New York-Heidelberg, 1979.

\item{[V1]} D. Voiculescu, {\eightsl The analogues of entropy and of
Fisher's information measure in free probability theory II,}
Invent.\ Math., {\eightbf 118} (1994), 411-440.

\item{[V2]} D. Voiculescu, {\eightsl The analogues of entropy and of
Fisher's information measure in free probability theory III: The
absence of Cartan subalgebras,} Geom.\ Funct.\ Anal.\ {\eightbf 6}
(1996) 172--199.

\item{[V3]} D. Voiculescu, {\eightsl Free entropy dimension $\le1$
for some generators of property T factors of type II$_1$},  J.
Reine Angew. Math. {\eightbf 514} (1999), 113--118.

\item{[VDN]} D. Voiculescu,
K. Dykema and A. Nica,  {\eightsl ``Free Random Variables,''} CRM
Monograph Series, vol. 1, AMS, Providence, R.I., 1992.
\end